\documentclass[]{article}

\usepackage[a4paper, total={6in, 8in}]{geometry}
\usepackage[]{algorithm2e}
\usepackage{amsmath}
\usepackage{amsfonts}
\usepackage{amsthm}
\usepackage{authblk}
\usepackage[ngerman, UKenglish]{babel}
\usepackage{bbm}
\usepackage{bigstrut}
\usepackage{caption}
\usepackage{cases}
\usepackage{comment}
\usepackage{float}
\usepackage[acronym]{glossaries}
\usepackage{graphicx, subfigure}
\usepackage[utf8]{inputenc}
\usepackage{mathtools}
\usepackage{moreverb}
\usepackage{multirow}
\usepackage{pgfplots}
\pgfplotsset{compat=1.3}
\usepackage{placeins}
\usepackage{tabularx}
\usepackage{tikz}

\makeglossaries
\DeclareMathAlphabet\mathbfcal{OMS}{cmsy}{b}{n}
\newacronym{cv}{CV}{cross-validation}
\newacronym{dc}{DC}{downward-closed}
\newacronym{hc}{HC}{hyperbolic cross}
\newacronym{lar}{LAR}{least angle regression}
\newacronym{lsr}{LSR}{least square regression}
\newacronym{mc}{MC}{Monte Carlo}
\newacronym{pce}{PCE}{polynomial chaos expansion}
\newacronym{pdf}{PDF}{probability density function}
\newacronym{psp}{PSP}{pseudo-spectral projection}
\newacronym{rv}{RV}{random variable}
\newacronym{rms}{RMS}{root-mean-square}
\newacronym[longplural={quantities of interest}]{qoi}{QoI}{quantity of interest}
\newacronym[longplural={sensitivity analyses}]{sa}{SA}{sensitivity analysis}
\newacronym{td}{TD}{total degree}
\newacronym{tp}{TP}{tensor product}
\newacronym{uq}{UQ}{uncertainty quantification}

\title{Adaptive Sparse Polynomial Chaos Expansions via Leja Interpolation}

\author[1,2,*]{Dimitrios Loukrezis}
\author[1,2]{Herbert De Gersem}

\affil[1]{{\small Institute for Accelerator Science and Electromagnetic Fields (TEMF), Technische Universit\"at Darmstadt \newline Schlossgartenstra\ss e 8, 64289 Darmstadt, Germany}}
\affil[2]{{\small Centre for Computational Engineering, Technische Universit\"at Darmstadt \newline Dolivostra\ss e 15, 64293 Darmstadt, Germany}}
\affil[*]{{\small Corresponding author (loukrezis@temf.tu-darmstadt.de)}}

\DeclareMathOperator*{\argmax}{argmax}
\date{}

\begin{document}

\maketitle

\providecommand{\keywords}[1]{\textbf{\textit{keywords--}} #1}

\begin{abstract}
This work suggests an interpolation-based stochastic collocation method for the non-intrusive and adaptive construction of sparse polynomial chaos expansions (PCEs).
Unlike pseudo-spectral projection and regression-based stochastic collocation methods, the proposed approach results in PCEs featuring one polynomial term per collocation point.
Moreover, the resulting PCEs are interpolating, i.e., they are exact on the interpolation nodes/collocation points.
Once available, an interpolating PCE can be used as an inexpensive surrogate model, or be post-processed for the purposes of uncertainty quantification and sensitivity analysis.
The main idea is conceptually simple and relies on the use of Leja sequence points as interpolation nodes.
Using Newton-like, hierarchical basis polynomials defined upon Leja sequences, a sparse-grid interpolation can be derived, the basis polynomials of which are unique in terms of their multivariate degrees.
A dimension-adaptive scheme can be employed for the construction of an anisotropic interpolation.
Due to the degree uniqueness, a one-to-one transform to orthogonal polynomials of the exact same degrees is possible and shall result in an interpolating PCE.
However, since each Leja node defines a unique Newton basis polynomial, an implicit one-to-one map between Leja nodes and orthogonal basis polynomials exists as well.
Therefore, the in-between steps of hierarchical interpolation and basis transform can be discarded altogether, and the interpolating PCE can be computed directly.
For directly computed, adaptive, anisotropic interpolating PCEs, the dimension-adaptive algorithm is modified accordingly.
A series of numerical experiments verify the suggested approach in both low and moderately high-dimensional settings, as well as for various input distributions.

\end{abstract}

\section{Introduction}
\label{sec:intro}
Studying a physical system often involves assessing the influence of uncertain parameters upon its behavior and outputs, the latter commonly called the \glspl{qoi}.
Typical study domains of this category are \gls{uq} \cite{smith2014, sullivan2015}, which aims to characterize  the \gls{qoi} statistically, and \gls{sa} \cite{saltelli2008}, the goal of which is to quantify the significance of each uncertain parameter regarding its impact on the \gls{qoi}. 
For both types of analysis, of particular importance to industrial scientists and engineers are non-intrusive methods, i.e., ones that allow the black-box use of the complex and often proprietary computer models which simulate the physical system under investigation.
Such models are typically computationally expensive, therefore, the cost of non-intrusive methods is usually measured by the number of model evaluations, accordingly, simulations, necessary until the desired accuracy has been reached.
Put otherwise, the less model evaluations needed for a sought accuracy, the more attractive the non-intrusive method. 
Several approaches that fall into this category are available, e.g., based on sampling \cite{caflisch1998, lemieux2009, loh1996, shields2015}, kriging \cite{dwight2009, lockwood2012, shimoyama2013}, or, more recently, neural networks \cite{tripathy2018, yang2018, zhu2019}.

If the functional relation between the model inputs and the \gls{qoi} is smooth, a spectral method \cite{lemaitre2010, xiu2010} can be used to derive a global polynomial approximation of it.
Once available, the polynomial approximation can be used as an inexpensive surrogate model and be sampled in the place of the original model, or be otherwise post-processed for \gls{uq} or \gls{sa} purposes \cite{babuska2010, buzzard2011, ghanem1991, nobile2008, sudret2008, xiu2005}.
For sufficiently smooth functions, a fast, even exponential, convergence in terms of approximation accuracy is to be expected.
The non-intrusive computation of spectral approximations is commonly accomplished with stochastic collocation methods \cite{eldred2009a, eldred2009, xiu2016}. 
Combined with sparse approximation techniques \cite{smolyak1963, temlyakov2015}, stochastic collocation methods may even break the curse of dimensionality \cite{chkifa2015a}, or at least enable studies with a comparably large number of parameters, typically up to moderately high dimensions \cite{blatman2010a, blatman2011, chkifa2014, conrad2013, doostan2011, hampton2015a, loukrezis2019, ma2009, migliorati2013a, nobile2008a, peng2014, winokur2016}.

One popular approach is to compute a (sparse) \gls{pce} \cite{ghanem1991, xiu2002}, the basis polynomials of which are orthogonal with respect to the joint \gls{pdf} which characterizes the model parameters. 
Orthogonal polynomials are chosen either according to the Wiener-Askey scheme \cite{xiu2002} or are numerically constructed \cite{soize2004, wan2006}, the latter option allowing the consideration of arbitrary input \glspl{pdf}.
Statistics and global sensitivity indices \cite{sobol2001} can be seamlessly computed by simply post-processing the \gls{pce}'s terms \cite{blatman2010a, ghanem1991, sudret2008}.
Non-intrusively, sparse \glspl{pce} are commonly computed with \gls{psp} \cite{conrad2013, constantine2012, winokur2016} or \gls{lsr} \cite{blatman2011, doostan2011, hampton2015a, migliorati2013a, peng2014} stochastic collocation methods.
Either to avoid aliasing errors \cite{conrad2013, constantine2012} or to establish stable regression problems \cite{chkifa2015, migliorati2013, migliorati2014}, both approaches result in less \gls{pce} terms than collocation points.
Equivalently, to construct an $M$-term \gls{pce}, $CM$, $C>1$, collocation points are needed.
Since each collocation point corresponds to a possibly costly model evaluation or simulation, it is of interest to keep their number to a minimum.

An $M$-term polynomial approximation can be constructed using exactly $M$ collocation points with an interpolation-based stochastic collocation method \cite{xiu2010, xiu2002}. 
However, with the exception of Lagrange interpolation schemes \cite{babuska2010, barthelmann2000, klimke2005, loukrezis2019, ma2009, nobile2008, nobile2008a, stoyanov2016, xiu2005}, the interpolation-based computation of spectral approximations is rarely seen in practice due to numerical stability issues \cite{trefethen2013}.
In the particular case of \glspl{pce}, the direct use of interpolation is not to be found in the literature.
Indirectly, a \gls{pce} can be derived by first computing a sparse-grid interpolation and then transforming the Lagrange polynomial basis to an orthogonal one \cite{buzzard2012, buzzard2013, porta2014}.
Nevertheless, such an approach requires multiple Lagrange polynomials of the same degree to be mapped to orthogonal ones of unique degrees.
Equivalently, multiple collocation points correspond to a single \gls{pce} basis polynomial. 
As a result, less \gls{pce} terms than collocation points are obtained, similar to the aforementioned \gls{psp} and \gls{lsr} collocation methods.

This work suggests an interpolation-based stochastic collocation method to construct sparse \glspl{pce} adaptively and using a single interpolation node, equivalently, collocation point, per polynomial term. 
The resulting \gls{pce} is an interpolating one, i.e., it is exact on the interpolation nodes.
The proposed method is based on Leja sequences \cite{leja1957}, which, additionally to other attractive properties, are nested and can be used for the construction of sparse interpolation grids \cite{narayan2014}.
The main idea is very simple and proceeds as follows.
By combining  Leja nodes with an interpolation scheme based on Newton-like, hierarchical polynomials \cite{chkifa2014}, a global approximation with basis polynomials of unique multivariate degrees can be constructed.
Additionally using a dimension-adaptive algorithm \cite{gerstner2003} results in an adaptively constructed, sparse-grid interpolation, while the use of weighted Leja nodes \cite{narayan2014} allows us to tailor the interpolation to any continuous input \gls{pdf}. 
In this case, a transform to an orthogonal polynomial basis has the distinct characteristic that the Newton basis polynomials are uniquely mapped to the orthogonal ones in terms of their degree, thus resulting in an interpolating \gls{pce}.
However, since each Leja node defines a unique Newton polynomial, it is implicitly mapped to a unique orthogonal polynomial as well. 
Therefore, the interpolation can be computed by directly using an orthogonal polynomial basis, i.e., without the in-between steps of a hierarchical interpolation and a basis transform.
The aforementioned dimension-adaptive algorithm can be accordingly modified as to employ orthogonal polynomials, while the adaptivity shall be guided by the \gls{pce} coefficients, which can be interpreted as sensitivity indicators.

The remainder of this paper is organized as follows.
Section~\ref{sec:prelim} presents the type of problem considered in this work and recalls briefly \glspl{pce} and interpolation-based stochastic collocation methods.
Section~\ref{sec:sai_pce} constitutes the main body of this paper and outlines the steps leading to adaptively constructed, sparse interpolating \glspl{pce}.
In Section~\ref{sec:numexp}, non-adaptive interpolating \glspl{pce} are compared against ones computed with \gls{psp} and \gls{lsr} methods, while adaptively constructed interpolating \glspl{pce} are verified in the context of moderately high-dimensional input uncertainty.
A conclusion on the findings of this work and suggestions on possible continuations are available in Section~\ref{sec:concl}.

\section{Preliminaries}
\label{sec:prelim}

\subsection{Problem Setting}
\label{subsec:prob_model}
We abstractly denote the functional relation between a model's input parameters and a scalar \gls{qoi} with $g\left(\mathbf{y}\right) \in \mathbb{R}$, where $\mathbf{y} \in \mathbb{R}^N$ is the parameter vector.
In this regard, function $g$ may refer to any kind of model, ranging from an analytical function to a complex numerical simulation.
In the context of the present work, $\mathbf{y} = \mathbf{Y}\left(\theta\right) \in \Xi$, $\theta \in \Theta$, is a realization of an $N$-dimensional \gls{rv} $\mathbf{Y} = \left(Y_1, \dots, Y_N\right)$, defined on the probability space $\left(\Theta, \Sigma, P\right)$ and characterized by the \gls{pdf} $\varrho: \Xi \rightarrow \mathbb{R}_{\geq 0}$, where $\Xi \subseteq \mathbb{R}^N$ denotes the image space, $\Theta$ the sample space, $\Sigma$ the set of events, and $P$ the probability measure.
Throughout this paper, we will assume that $\mathbf{Y}$ consists of independent \glspl{rv}, such that $\varrho(\mathbf{y}) = \prod_{n=1}^N \varrho_n(y_n)$, $\Theta = \Theta_1 \times \cdots \times \Theta_N$, and $\Xi = \Xi_1 \times \cdots \times \Xi_N$.
Nevertheless, dependencies can also be handled with suitable transformations \cite{feinberg2018, jakeman2019, lebrun2009}. 

Under the assumption that $g$ is a smooth function, our goal is to compute a spectral approximation 
\begin{equation}
\label{eq:spectral_approx}
g(\mathbf{y}) \approx \widetilde{g}(\mathbf{y}) = \sum_{m=1}^M c_m \Psi_m(\mathbf{y}),
\end{equation}
where $c_m$ are scalar coefficients and $\Psi_m$ multivariate global polynomials, i.e., ones defined over the entire image space $\Xi$.
Given sufficient approximation accuracy, the polynomial model \eqref{eq:spectral_approx} may be used as a surrogate to the original model and replace it in computationally demanding tasks to reduce the overall cost.
Moreover, depending on the polynomial basis $\left\{\Psi_m\right\}_{m=1}^M$, statistical measures and sensitivity indices can be inexpensively estimated by post-processing the approximation's coefficients and polynomials \cite{babuska2010, buzzard2011, ghanem1991, nobile2008, sudret2008, xiu2005}. 

\subsection{Polynomial Chaos Expansions}
\label{subsec:pce}
The polynomial approximation \eqref{eq:spectral_approx} is a \gls{pce} if the polynomials $\Psi_m$, $m=1,\dots,M$, are orthogonal with respect to the joint \gls{pdf}, i.e., if they satisfy the property
\begin{equation}
\label{eq:orthNd}
\mathbb{E}\left[\Psi_m \Psi_n\right] = \int_\Xi \Psi_m\left(\mathbf{y}\right) \Psi_n\left(\mathbf{y}\right) \varrho\left(\mathbf{y}\right) \mathrm{d}\mathbf{y} = \gamma_m \delta_{mn},
\end{equation}
where $\delta_{mn}$ is the Kronecker delta and $\gamma_m$ a normalization factor.
The orthogonal polynomials are commonly chosen according to the Wiener-Askey scheme \cite{xiu2002} or are numerically constructed \cite{soize2004, wan2006}.
Under the independence assumption, the multivariate orthogonal polynomials are given as  
\begin{equation}
\label{eq:polyNd}
\Psi_m \equiv \Psi_\mathbf{p}(\mathbf{y}) = \prod_{n=1}^N \psi_n^{p_n} (y_n),
\end{equation}
where $\psi_n^{p_n}$ are univariate polynomials of degree $p_n \in \mathbb{Z}_{\geq 0}$, which are orthogonal with respect to the corresponding univariate \gls{pdf}, i.e., they satisfy the property
\begin{equation}
\label{eq:orth1d}
\mathbb{E}\left[\psi_n^{p_n}\psi_n^{q_n}\right] = \int_{\Xi} \psi_n^{p_n}(y_n) \psi_n^{q_n}(y_n) \varrho_n(y_n) \, \mathrm{d}y_n = \gamma_{p_n} \delta_{p_nq_n}. 
\end{equation}
Accordingly, $\mathbf{p} = \left(p_1, \dots, p_N\right)$ is a multi-index equivalent to the multivariate polynomial degree.
We note that a single index $m$ is uniquely associated to a multi-index $\mathbf{p}$, e.g., as in \cite{babuska2010}.

Using the multi-index notation, formula \eqref{eq:spectral_approx} can be equivalently written as
\begin{equation}
\label{eq:spectral_approx_multi_index}
g(\mathbf{y}) \approx \widetilde{g}(\mathbf{y}) = \sum_{\mathbf{p} \in \Lambda} c_{\mathbf{p}} \Psi_{\mathbf{p}}(\mathbf{y}),
\end{equation}
where $\Lambda$ is a multi-index set with cardinality $\#\Lambda = M$.
\Gls{tp}, \gls{td}, and \gls{hc} multi-index sets are commonly encountered in the literature, see, e.g., \cite{babuska2010, blatman2011, constantine2012}, all of which fall into the category of \gls{dc} sets \cite{gerstner2003}, i.e., they satisfy the property
\begin{equation}
\label{eq:monotonicity}
\forall \mathbf{p} \in \Lambda \Rightarrow \mathbf{p}-\mathbf{e}_n \in \Lambda, \forall n = 1,2,\dots,N, \: \text{with} \: p_n > 0,
\end{equation}
where $\mathbf{e}_n$ is the $n$-th unit vector.
Other than \gls{tp} ones, \gls{dc} multi-index sets are often employed in the context of sparse \glspl{pce} \cite{conrad2013, constantine2012, migliorati2013a, winokur2016}.
Arbitrarily shaped sets can also be used \cite{blatman2010, blatman2010a}.

The post-processing of \glspl{pce} has been explained thoroughly elsewhere, see, e.g., \cite{blatman2010a, ghanem1991, sudret2008}. 
We remind briefly some of its aspects, with an emphasis on \gls{pce}-based Sobol \gls{sa} \cite{saltelli2008, sobol2001}, which is closely connected to the dimension-adaptive construction of interpolating \glspl{pce} discussed in Section~\ref{subsec:basis_equiv}.
In particular, the first \gls{pce} coefficient, i.e., the one corresponding to the multi-index $\mathbf{0} = \left(0, \dots, 0\right)$, provides an estimate for the expected value of the \gls{qoi}, such that
\begin{equation}
\label{eq:pce_mean}
\mathbb{E}\left[g\right] = c_{\mathbf{0}}.
\end{equation}
Assuming that the \gls{pce} basis is orthonormal, i.e., $\gamma_m=1$, $m=1,\dots,M$, in formula \eqref{eq:orthNd}, the variance of the \gls{qoi} is estimated by
\begin{equation}
\mathbb{V}\left[g\right] = \sum_{\mathbf{p} \in \Lambda \setminus \mathbf{0}} c_{\mathbf{p}}^2.
\end{equation}
Regarding the variance-based Sobol \gls{sa}, we focus here on the so-called first-order and total-order Sobol indices. 
The former quantify the impact of a \gls{rv} upon the \gls{qoi} assuming no interaction with the other \glspl{rv}, which are kept constant.
In the latter case, the \gls{rv} is assumed to interact with all other \glspl{rv}.
Defining the first-order and total-order multi-index sets with respect to the $n$-th \gls{rv} as
\begin{align}
\Lambda_n^{\text{f}} &= \{\mathbf{p} \in \Lambda \; : \; p_n \neq 0 \:\: \text{and} \:\:p_m = 0, m \neq n\}, \\
\Lambda_n^{\text{t}} &= \{\mathbf{p} \in \Lambda \; : \; p_n \neq 0\},
\end{align}
we compute the corresponding partial variances
\begin{align}
\mathbb{V}_n^{\text{f}}\left[g\right] = \sum_{\mathbf{p}  \in \Lambda_n^{\text{f}}} c_{\mathbf{p}}^2, \hspace{1em}
\mathbb{V}_n^{\text{t}}\left[g\right] = \sum_{\mathbf{p}  \in \Lambda_n^{\text{t}}} c_{\mathbf{p}}^2,
\end{align}
and, finally, the first- and total-order Sobol indices
\begin{align}
S_n^{\text{f}} = \frac{\mathbb{V}_n^{\text{f}}\left[g\right]}{\mathbb{V}\left[g\right]}, \hspace{1em}
S_n^{\text{t}} = \frac{\mathbb{V}_n^{\text{t}}\left[g\right]}{\mathbb{V}\left[g\right]}.
\end{align}
Sobol indices corresponding to in-between interaction orders are here omitted, but can be computed in a similar way.

\subsection{Interpolation-Based Stochastic Collocation}
\label{subsec:interpolation}
One possible way to compute \eqref{eq:spectral_approx} non-intrusively, is to use an interpolation-based stochastic collocation approach \cite{xiu2010, xiu2016}. 
Given an interpolation grid, i.e., a set of prescribed nodes $\left\{ \mathbf{y}_m \right\}_{m=1}^M$ along with the corresponding model evaluations $\left\{ g(\mathbf{y}_m) \right\}_{l=1}^M$, the interpolation condition $g(\mathbf{y}_m) = \widetilde{g}(\mathbf{y}_m)$, $m=1, \dots, M$, leads to the linear system of equations
\begin{equation}
\label{eq:interpo_system_eq}
\mathbf{A} \mathbf{c} = \mathbf{g},
\end{equation}
where $\mathbf{A} \in \mathbb{R}^{M \times M}$ with elements $a_{mn} = \Psi_{n}\left(\mathbf{y}_m\right)$, $\mathbf{g} = \left(g(\mathbf{y}_1), \dots, g\left(\mathbf{y}_M\right)\right)^\top \in \mathbb{R}^M$, and $\mathbf{c} = \left(c_1, \dots, c_M\right)^\top \in \mathbb{R}^M$.
Solving \eqref{eq:interpo_system_eq} for $\mathbf{c}$ yields the coefficients of the spectral approximation \eqref{eq:spectral_approx}.

In principle, the interpolating polynomial is unique for a given interpolation grid and can be represented using any polynomial basis \cite{gander2005}. 
E.g., using an orthogonal basis with respect to the \gls{pdf}, an interpolating \gls{pce} would be obtained. 
In practice, however, due to numerical stability issues \cite{trefethen2013}, interpolation-based stochastic collocation methods are avoided. 
One exception is the so-called Lagrange interpolation approach, where, most commonly, the global polynomial approximation \eqref{eq:spectral_approx} is given as a linear combination of \gls{tp} univariate Lagrange interpolation rules.
The method is known as sparse-grid interpolation or as the sparse-grid stochastic collocation method, and is briefly introduced next. 
Extensive presentations are available in a number of works, see, e.g., \cite{babuska2010, barthelmann2000, nobile2008, xiu2005}. 

With respect to the $n$-th \gls{rv}, a univariate Lagrange interpolation rule is defined by an interpolation level $i_n \in \mathbb{Z}_{\geq 0}$ and the corresponding interpolation grid $G_{n, i_n} = \left\{y_{n, i_n}^{(j_n)}\right\}_{j_n=1}^{m_n(i_n)}$, where the monotonically increasing level-to-nodes function $m_n(i_n)$, $m_n(0)=1$, relates the interpolation level to the size of the grid.
The univariate Lagrange interpolation reads
\begin{equation}
\label{eq:lagr_interp_1d}
\mathcal{I}_{n, i_n} \left[g\right]\left(y_n\right) = \sum_{j_n=1}^{m_n(i_n)} g\left(y_{n, i_n}^{(j_n)}\right) l_{n, i_n}^{(j_n)}\left(y\right),
\end{equation}
where $l_{n, i_n}^{(j_n)}$ are univariate Lagrange polynomials of degree $m_n(i_n)-1$, defined as
\begin{equation}
\label{eq:lagr_poly_1d}
l_{n, i_n}^{(j_n)}\left(y_n\right) = 
\prod_{k=1, k \neq j_n}^{m_n(i_n)} \frac{y_n - y_{n, i_n}^{(k)}}{y_{n,i_n}^{(j_n)} - y_{n, i_n}^{(k)}}.
\end{equation}
It is generally desirable to use nested grids, such that $G_{n, i_n-1} \subset G_{n, i_n}$, $i_n > 0$.
Then, in the case of an interpolation level refinement from $i_n$ to $i_n+1$, the model evaluations upon the nodes of $G_{i_n}$ are readily available.
Moreover, nested univariate grids are crucial to the construction of sparse grids in the multivariate case \cite{barthelmann2000}.

For the multivariate sparse-grid interpolation, we employ a multi-index set $\Lambda$ holding multi-indices $\mathbf{i} = \left(i_1, \dots, i_N\right)$ with the corresponding univariate interpolation levels.
The multi-index set $\Lambda$ is typically \gls{dc} \cite{gerstner2003}.
Defining the univariate hierarchical surplus operator as
\begin{equation}
\label{eq:hs_operator_1d}
\Delta_{n, i_n} = \mathcal{I}_{n, i_n} - \mathcal{I}_{n, i_n - 1},
\end{equation}
where $\mathcal{I}_{n,-1}[g](y_n)=0$, and its \gls{tp} counterpart as 
\begin{equation}
\label{eq:hs_operator_nd}
\Delta_{\mathbf{i}} = \Delta_{1, i_1} \otimes \cdots \otimes \Delta_{N, i_N}, 
\end{equation}
where $\otimes$ denotes a tensor product, the sparse-grid interpolation reads
\begin{equation}
\label{eq:lagr_interp_nd}
\mathcal{I}_{\Lambda}[g](\mathbf{y}) = \sum_{\mathbf{i} \in \Lambda} \Delta_{\mathbf{i}}[g](\mathbf{y}).
\end{equation}
The corresponding sparse grid is given as a union of tensor grids, such that
\begin{equation}
\label{eq:interp_grid_nd}
G_{\Lambda} = \bigcup_{\mathbf{i} \in \Lambda} G_{\mathbf{i}} = \bigcup_{\mathbf{i} \in \Lambda} \left(G_{1,i_1} \times \cdots \times G_{N, i_N}\right).
\end{equation}
Depending on the form of the multi-index set $\Lambda$, sparse grids are either isotropic \cite{babuska2010, barthelmann2000, nobile2008, xiu2005} or anisotropic \cite{chkifa2014, klimke2005, loukrezis2019, ma2009, narayan2014, nobile2008a, schillings2013, stoyanov2016}, the latter being especially advantageous in the case of parameters with an unequal influence upon the \gls{qoi}.

Formula \eqref{eq:lagr_interp_nd} can be equivalently written in the form of \eqref{eq:spectral_approx}, where $\Psi_m$ are now multivariate Lagrange polynomials
\begin{equation}
\label{eq:lagr_poly_nd}
L_{\mathbf{i}}^{(\mathbf{j})} = \prod_{n=1}^N l_{n, i_n}^{(j_n)},
\end{equation}
of multivariate degrees $\left(m_1(i_1)-1, \dots, m_N(i_N)-1\right)$.
Contrary to formula \eqref{eq:spectral_approx_multi_index}, it now holds that $\#G_{\Lambda} = M$ instead of $\#\Lambda=M$, since each interpolation node defines a single polynomial and a multi-index may correspond to more than one interpolation nodes.
We also note that formula \eqref{eq:lagr_interp_nd} does not necessarily satisfy the interpolation property
\begin{equation}
\label{eq:interp_prop_nd}
\mathcal{I}_{\Lambda}[g]\left(\mathbf{y}_{\mathbf{i}}^{(\mathbf{j})}\right) = g\left( \mathbf{y}_{\mathbf{i}}^{(\mathbf{j})} \right), \hspace{1em} \mathbf{y}_{\mathbf{i}}^{(\mathbf{j})} = \left(y_{1,i_1}^{(j_1)}, \dots, y_{N,i_N}^{(j_N)}\right) \in G_\Lambda,
\end{equation}
unless the underlying univariate interpolation rules \eqref{eq:lagr_interp_1d} are based on nested grids \cite{barthelmann2000}.

\section{Adaptive Sparse Polynomial Chaos Interpolation}
\label{sec:sai_pce}

\subsection{Leja Sequences}
\label{subsec:leja_seq_seq}
Leja sequences constitute the main component of the suggested approach for computing an interpolating \gls{pce}.
In its original form \cite{leja1957}, a Leja node sequence  $\left\{y^{(j)}\right\}_{j \geq 0}$, $y^{(j)} \in \left[-1,1\right]$, is obtained by solving the optimization problem
\begin{equation}
\label{eq:lejaopt_unweighted}
 y^{(j)} = \argmax_{y \in \left[-1, 1\right]} \prod_{k=0}^{j-1}\left|y - y^{(k)}\right|,
\end{equation}
for each node $y^{(j)}$, $j \geq 1$.
The initial node $y^{(0)} \in \left[-1,1\right]$ can be chosen arbitrarily, therefore Leja sequences are not unique.
The weighted counterpart of \eqref{eq:lejaopt_unweighted} uses a univariate \gls{pdf} $\varrho\left(y\right)$ as a weight function and transforms the optimization problem to 
\begin{equation}
\label{eq:lejaopt_weighted}
 y^{(j)} = \argmax_{y \in \Xi} \sqrt{\varrho\left(y\right)} \prod_{k=0}^{j-1}\left|y - y^{(k)}\right|,
\end{equation}
thus yielding a so-called weighted Leja sequence \cite{narayan2014}.
As in the unweighted case, the initial node can be chosen arbitrarily in the image space $\Xi$.

Leja nodes have a number of desirable properties.
We focus here on the properties that are relevant in the context of this work. 
A more extensive discussion can be found in \cite{narayan2014}.
First, Leja sequences comprise suitable interpolation and quadrature nodes, as indicated by tests with uniform \cite{chkifa2014, loukrezisPhD, loukrezis2019, narayan2014, nobile2015, schillings2013}, normal \cite{narayan2014}, and more general distributions \cite{farcas2019, georg2018, loukrezisPhD, loukrezis2019, vandenbos2018}.
Additionally, due to their definition in \eqref{eq:lejaopt_unweighted} and \eqref{eq:lejaopt_weighted}, Leja sequences are nested.
Due to the nestedness property, Leja nodes can be employed in the - possibly adaptive - construction of sparse interpolation and quadrature grids \cite{chkifa2014, georg2018, loukrezisPhD, loukrezis2019, narayan2014, nobile2015, schillings2013}.

\subsection{Hierarchical Leja Interpolation}
\label{subsec:hierarch_leja_interp}
In connection to the Lagrange interpolation method presented in Section~\ref{subsec:interpolation}, we now assume that the level-$i_n$ univariate interpolation grid coincides with a Leja sequence given as in \eqref{eq:lejaopt_unweighted} or \eqref{eq:lejaopt_weighted}, such that $G_{i_n} = \left\{ y_{n}^{(j_n)} \right\}_{j_n=0}^{i_n}$.
In that case, the corresponding level-to-nodes function is $m_n(i_n) = i_n + 1$. 
Moreover, the interpolation grids are nested and differ by a single node, i.e., $G_{i_n-1} \subset G_{i_n}$ with $G_{i_n} = G_{i_n-1} \cup y_n^{(i_n)}$, $i_n > 0$.

Next, we abandon the classical Lagrange interpolation \eqref{eq:lagr_interp_1d} in favor of a hierarchical/sequential one \cite{chkifa2014}.
To that end, we introduce the modified Newton polynomials
\begin{equation}
\label{eq:hier_poly_1d}
\nu_{n}^{i_n}\left(y_n\right) = 
\begin{dcases}
\prod_{k=0}^{i_n-1} \frac{y_n-y_{n}^{(k)}}{y_{n}^{(i_n)} - y_{n}^{(k)}},  &i_n \neq 0, \\
1,  &i_n=0.
\end{dcases}
\end{equation}
Essentially, a modified Newton polynomial $\nu_n^{i_n}$ coincides with the Lagrange polynomial $l_{n, i_n}^{(i_n+1)}$ defined for the same Leja sequence.
However, the degree of $\nu_{n}^{i_n}$ is unique and equal to $i_n$, such that no two polynomials defined for the same Leja sequence can be of equal degree.
Equivalently, each Leja node defines and corresponds to a single polynomial.
Moreover, Newton polynomials do not change if additional nodes are added to the interpolation grid.
One can then derive the hierarchical/sequential univariate interpolation
\begin{subequations}
\label{eq:hier_interp_1d}	
\begin{align}
\mathcal{I}_{n, i_n}\left[g\right]\left(y_n\right)  
&= \sum_{j_n=0}^{i_n} s_{j_n} \nu_{n}^{j_n}(y_n) \\
&= \sum_{j_n=0}^{i_n} \Delta_{n, j_n}\left[g\right](y_n) \\
&=  \mathcal{I}_{n, i_n-1}\left[g\right]\left(y_n\right) + \left( g\left(y_n^{(i_n)}\right) - \mathcal{I}_{n, i_n-1}\left[g\right]\left(y_n^{(i_n)}\right)\right) \nu_{n}^{i_n}(y_n),
\end{align}
\end{subequations}
where the coefficients $s_{j_n}$ are known as the hierarchical surpluses. 

Moving to the multivariate case, each multi-index $\mathbf{i}$ defines now a single multi-dimensional Leja node  $\mathbf{y}^{(\mathbf{i})} = \left(y_{1}^{(i_1)}, \dots, y_{N}^{(i_N)}\right)$ and a single multivariate modified Newton polynomial 
\begin{equation}
\label{eq:hierarch_Nd}
N_{\mathbf{i}}\left(\mathbf{y}\right) = \prod_{n=1}^N \nu_{n}^{i_n}(y_n),
\end{equation}
the multivariate degree of which is unique and equal to $\mathbf{i}$, similarly to the univariate case.
Accordingly, employing the sparse-grid interpolation formula \eqref{eq:lagr_interp_nd}, each hierarchical surplus operator $\Delta_{\mathbf{i}}$ shall add a single polynomial term to the interpolation and a single Leja node to the sparse grid. 
It therefore holds that $\#\Lambda = \#G_{\Lambda}$ and, similarly to formula \eqref{eq:spectral_approx_multi_index}, $\#\Lambda = M$.
Consequently, the summation in \eqref{eq:lagr_interp_nd} is performed over a sequence of \gls{dc} multi-index sets $\left(\Lambda_k\right)_{k=1}^{\#\Lambda}$, such that $\#\Lambda_k = \#G_{\Lambda_k} = k$, $\Lambda_{k} = \Lambda_{k-1} \cup \mathbf{i}_k$, and $G_{\Lambda_k} = G_{\Lambda_{k-1}} \cup \mathbf{y}^{(\mathbf{i}_k)}$, where $\Lambda_0 = G_{\Lambda_0} = \emptyset$. 
Then, the interpolation with respect to each multi-index set $\Lambda_k$ can be written in the hierarchical/sequential form
\begin{subequations}
\label{eq:hier_interp_nd}
\begin{align}
\mathcal{I}_{\Lambda_k}[g](\mathbf{y}) &= \sum_{\mathbf{i} \in \Lambda_k} s_{\mathbf{i}} N_{\mathbf{i}}(\mathbf{y}) \\
&= \sum_{\mathbf{i} \in \Lambda_k} \Delta_{\mathbf{i}}[g](\mathbf{y}) \\
& = \mathcal{I}_{\Lambda_{k-1}}\left[g\right]\left(\mathbf{y}\right) + \left(g\left(\mathbf{y}^{(\mathbf{i}_k)}\right) - \mathcal{I}_{\Lambda_{k-1}}\left[g\right]\left(\mathbf{y}^{(\mathbf{i}_k)}\right)\right) N_{\mathbf{i}_k}\left(\mathbf{y}\right).
\end{align}
\end{subequations}

In principle, the sequence $\left(\Lambda_k\right)_{k \geq 1}$ with $\#\Lambda_k = k$ can be constructed by appropriately ordering the multi-indices of any \gls{dc} multi-index set, such as a \gls{td} or a \gls{hc} one.
Alternatively, a greedy-adaptive algorithm can be employed, similar to the dimension-adaptive approach originally presented in \cite{gerstner2003} for quadrature purposes, in order to construct an anisotropic sparse-grid interpolation. 
We tailor the idea of dimension-adaptivity to the special case of Leja-sequence-based interpolation grids, as depicted in Algorithm~\ref{algo:dimadapt_hierarchical}.

\begin{algorithm}[t!]
	\SetAlgoLined
	\KwData{function $g\left(\mathbf{y}\right)$, initial \gls{dc} multi-index set $\Lambda^{\text{init}}$, tolerance $\epsilon$, simulation budget $B$.}
	\KwResult{\gls{dc} multi-index set $\Lambda$, sparse grid $G_{\Lambda}$, sparse-grid interpolation $\mathcal{I}_{\Lambda}\left[g\right]$.}
	$k = \#\Lambda^{\text{init}}$. \\ $\Lambda_k \leftarrow \Lambda^{\text{init}}$.\\
	\While{\textsc{True}}{
		$\Lambda_k^{\text{adm}} = \left\{\mathbf{i} : \mathbf{i} \not \in \Lambda_k \:\: \text{and} \:\: \Lambda_k \cup \mathbf{i} \:\: \text{is DC}\right\}$. \\
		$ s_{\mathbf{i}} = g\left(\mathbf{y}^{(\mathbf{i})}\right) - \mathcal{I}_{\Lambda_k}\left[g\right]\left(\mathbf{y}^{(\mathbf{i})}\right)$, $\forall \mathbf{i} \in \Lambda_k^{\text{adm}}$. \\
		\If{$\#\Lambda_k + \#\Lambda_k^{\mathrm{adm}} \geq B$ \textsc{or} $\sum_{\mathbf{i} \in \Lambda_k^{\mathrm{adm}}} \left|s_{\mathbf{i}}\right| \leq \epsilon$}{Exit while-loop.}
		$ \mathbf{i}^*  = \argmax_{\mathbf{i} \in \Lambda_k^{\text{adm}}}\left|s_{\mathbf{i}}\right|$. \\
		$\Lambda^* \leftarrow \Lambda_k \cup \mathbf{i}^*$. \\ 
		$k \leftarrow k+1$. \\ $\Lambda_k \leftarrow \Lambda^*$.
	}
	$\Lambda \leftarrow \Lambda_k \cup \Lambda_k^{\text{adm}}$.
	\caption{Dimension-adaptive hierarchical Leja interpolation.} 
	\label{algo:dimadapt_hierarchical}
\end{algorithm}

The algorithm is typically initialized with $\Lambda_1 = \left\{ \mathbf{i}_1 = \mathbf{0} \right\}$, however, any initial \gls{dc} set can be used.
At each step, a \gls{dc} set $\Lambda_k$ with $\#\Lambda_k = k$ is available, along with the Leja-based sparse grid $G_{\Lambda_k}$ and the interpolation $\mathcal{I}_{\Lambda_k}\left[g\right](\mathbf{y})$.
The corresponding admissible set, i.e., the set of multi-indices which satisfy property \eqref{eq:monotonicity} if added to $\Lambda_k$, is defined as 
\begin{equation}
\label{eq:admset}
\Lambda_k^{\text{adm}} = \left\{\mathbf{i} : \mathbf{i} \not \in \Lambda_k \:\: \text{and} \:\: \Lambda_k \cup \mathbf{i} \:\: \text{is DC}\right\}.
\end{equation}
Each multi-index $\mathbf{i} \in \Lambda_k^{\text{adm}}$ is uniquely associated to a single admissible Leja node $\mathbf{y}^{(\mathbf{i})}$.
Similar to \cite{gerstner2003}, we connect the contribution of an admissible multi-index to the magnitude of the corresponding hierarchical surplus $s_{\mathbf{i}}$, which can be computed as in \eqref{eq:hier_interp_nd}. 
Therefore, the multi-index set $\Lambda_k$ shall be enriched with the multi-index
\begin{equation}
\label{eq:next_multi_index}
\mathbf{i}^* = \argmax_{\mathbf{i} \in \Lambda_k^{\text{adm}}} \left|s_{\mathbf{i}}\right|,
\end{equation}
i.e., the one corresponding to the maximum contribution.

At any given step of the algorithm, the total number of model evaluations is equal to $\#\Lambda_k + \#\Lambda_k^{\text{adm}}$. 
The procedure is continued iteratively until a simulation budget $B$ is reached, or until the total contribution of the admissible set is below a tolerance $\epsilon$.
After the termination of the iterations, the final interpolation is constructed using all multi-indices in $\Lambda_k \cup \Lambda_k^{\text{adm}}$, since the hierarchical surpluses corresponding to the admissible nodes have already been computed.

\subsection{Basis Equivalence and Interpolating PCEs}
\label{subsec:basis_equiv}
Given a hierarchical interpolation in the form of \eqref{eq:hier_interp_nd}, we may use the multi-index set $\Lambda$ to construct the orthonormal polynomial basis $\left\{\Psi_{\mathbf{i}}\right\}_{\mathbf{i} \in \Lambda}$. 
For both types of basis polynomials, the respective multivariate polynomial degrees coincide with the multi-indices $\mathbf{i} \in \Lambda$. 
Hence, there exists a one-to-one relation between the modified Newton and the orthogonal polynomials in terms of their degrees.
We therefore obtain two equivalent representations of the interpolating polynomial upon the sparse grid $G_\Lambda$, i.e.,
\begin{equation}
\label{eq:interp_equivalence}
\mathcal{I}_\Lambda\left[g\right]\left(\mathbf{y}\right) = \sum_{\mathbf{i} \in \Lambda} s_{\mathbf{i}} N_{\mathbf{i}}(\mathbf{y}) = \sum_{\mathbf{i} \in \Lambda} c_{\mathbf{i}} \Psi_{\mathbf{i}}(\mathbf{y}),
\end{equation}
where the second representation is an interpolating \gls{pce}.
The unknown \gls{pce} coefficients can be computed simply by enforcing the interpolation property and solving the system \eqref{eq:interpo_system_eq}, where the elements of the right-hand-side vector $\mathbf{g}$ correspond to the already available model evaluations upon the nodes of the Leja sparse grid.

Essentially, the procedure described above is a basis transform, similar in spirit to the ones employed in \cite{buzzard2012, buzzard2013, porta2014} for the case of uniformly distributed inputs, Clenshaw-Curtis nodes, and Lagrange interpolation.
In those works, a mapping from the multiple Lagrange polynomials with an equal degree to the unique orthogonal ones is needed.
Accordingly, multiple collocation points correspond to a single orthogonal polynomial.
On the contrary, here, the aforementioned one-to-one map between Newton and orthogonal basis polynomials results in a much simplified basis transform, while, more importantly, the number of \gls{pce} terms is exactly equal to the number of interpolation nodes.
Moreover, the \gls{pce} is exact on the collocation points, i.e., it satisfies the interpolation property.
Additionally, the use of weighted Leja nodes allows us to consider any continuous input \gls{pdf} when constructing the initial interpolation.

Nevertheless, the interpolating \gls{pce} can be computed without the in-between steps of first deriving a hierarchical interpolation and then performing a basis transform.
As shown in Section~\ref{subsec:hierarch_leja_interp}, each Leja node $\mathbf{y}^{(\mathbf{i})}$, $\mathbf{i} \in \Lambda$, defines a unique modified Newton polynomial $N_{\mathbf{i}}$.
Therefore, based on equation \eqref{eq:interp_equivalence} and the aforementioned equivalence between the two polynomial bases, each Leja node is implicitly mapped to a unique orthogonal basis polynomial as well. 
Consequently, given a multi-index set $\Lambda$, the \gls{pce} can be computed by solving the linear system \eqref{eq:interpo_system_eq}, where the system matrix $\mathbf{A}$ is formed based on the orthogonal polynomial basis $\left\{\Psi_{\mathbf{i}}\right\}_{\mathbf{i} \in \Lambda}$ and the Leja points $\mathbf{y}^{(\mathbf{i})}$, $\mathbf{i} \in \Lambda$, while the right-hand-side vector $\mathbf{g}$ is obtained by evaluating the model on the Leja points.

The interpolating \gls{pce} can also be constructed adaptively, by suitably modifying Algorithm~\ref{algo:dimadapt_hierarchical}.
According to the \gls{pce}-based \gls{sa} discussed in Section~\ref{subsec:pce}, for a  normalized polynomial basis $\left\{\Psi_{\mathbf{i}}\right\}_{\mathbf{i} \in \Lambda}$, the value $c_{\mathbf{i}}^2$ corresponds to the partial variance due to the multi-index $\mathbf{i}$.
Therefore, each coefficient $c_{\mathbf{i}}$ can be interpreted as a sensitivity indicator regarding the corresponding multi-index and replace the hierarchical surpluses $s_{\mathbf{i}}$ in the criterion \eqref{eq:next_multi_index}, which is used for the expansion of the multi-index set.
The modified procedure is presented in Algorithm~\ref{algo:dimadapt_orthogonal}, where orthogonal polynomials are used instead of modified Newton ones.
Since the orthogonal basis is not a hierarchical one, instead of computing hierarchical surpluses, the system 
\begin{equation}
 \label{eq:interpo_system_algo}
 \mathbf{A}_{\text{sys}} \mathbf{c}_{\text{sys}} = \mathbf{g}_{\text{sys}},
\end{equation} 
must be solved at each step of the algorithm, where $\mathbf{A}_{\text{sys}} \in \mathbb{R}^{K \times K}$, $\mathbf{g}_{\text{sys}} \in \mathbb{R}^K$, and $\mathbf{c}_{\text{sys}} \in \mathbb{R}^K$, $K = \#\left(\Lambda_k \cup \Lambda_k^{\mathrm{adm}}\right)$, are defined similarly to the linear system \eqref{eq:interpo_system_eq}, for all multi-indices in the set $\Lambda_k \cup \Lambda_k^{\mathrm{adm}}$.
The multi-index set $\Lambda_k $ is then expanded with the admissible multi-index corresponding to the maximum sensitivity, equivalently, to the maximum value $\left|c_{\mathbf{i}}\right|$, $\mathbf{i} \in \Lambda_k^{\mathrm{adm}}$.

Compared to Algorithm~\ref{algo:dimadapt_hierarchical}, the modified algorithm comes at the expense of solving a linear system at each step of the algorithm.
Nevertheless, for a reasonable number of polynomials, equivalently, model evaluations, the cost overhead is not significant.
E.g., in most practical applications, the simulation budget is typically in the order of hundreds or few thousands.
Especially in the case of computationally expensive simulations, the cost due to solving the linear system can essentially be regarded as negligible.
\begin{algorithm}[t!]
	\SetAlgoLined
	\KwData{function $g\left(\mathbf{y}\right)$, initial \gls{dc} multi-index set $\Lambda^{\text{init}}$, tolerance $\epsilon$, simulation budget $B$.}
	\KwResult{\gls{dc} multi-index set $\Lambda$, sparse grid $G_{\Lambda}$, sparse-grid interpolation $\mathcal{I}_{\Lambda}\left[g\right]$.}
	$k = \#\Lambda^{\text{init}}$. \\ $\Lambda_k \leftarrow \Lambda^{\text{init}}$.\\
	\While{\textsc{True}}{
		$\Lambda_k^{\text{adm}} = \left\{\mathbf{i} : \mathbf{i} \not \in \Lambda_k \:\: \text{and} \:\: \Lambda_k \cup \mathbf{i} \:\: \text{is DC}\right\}$. \\
		$\Lambda_{\mathrm{sys}} = \Lambda_k \cup \Lambda_k^{\text{adm}}$. \\
		$\mathbf{c}_{\mathrm{sys}} = \mathbf{A}_{\mathrm{sys}}^{-1} \mathbf{g}_{\mathrm{sys}}$. \\
		\If{$\#\Lambda_k + \#\Lambda_k^{\mathrm{adm}} \geq B$ \textsc{or} $\sum_{\mathbf{i} \in \Lambda_k^{\mathrm{adm}}} \left|c_{\mathbf{i}}\right| \leq \epsilon$}{Exit while-loop.}
		$ \mathbf{i}^*  = \argmax_{\mathbf{i} \in \Lambda_k^{\text{adm}}}\left|c_{\mathbf{i}}\right|$. \\
		$\Lambda^* \leftarrow \Lambda_k \cup \mathbf{i}^*$. \\ 
		$k \leftarrow k+1$. \\ $\Lambda_k \leftarrow \Lambda^*$.
	}
	$\Lambda \leftarrow \Lambda_k \cup \Lambda_k^{\text{adm}}$.
	\caption{Dimension-adaptive, Leja-based, polynomial chaos interpolation.} 
	\label{algo:dimadapt_orthogonal}
\end{algorithm}

\section{Numerical Experiments}
\label{sec:numexp}
The models used in the numerical experiments of this section are freely accessible in a dedicated online repository\footnote{\texttt{https://github.com/dlouk/UQ\_benchmark\_models}}.
Algorithms~\ref{algo:dimadapt_hierarchical} and \ref{algo:dimadapt_orthogonal} have been implemented in the in-house DALI (Dimension-Adaptive Leja Interpolation) software\footnote{\texttt{https://github.com/dlouk/DALI3}} \cite{loukrezis2019}.
The Leja nodes, either weighted or unweighted, are provided by the Chaospy Python library \cite{feinberg2015}.
For the orthogonal polynomials, we use the OpenTURNS C++/Python library \cite{baudin2017}.
All numerical experiments have been carried out on a desktop computer equipped with  an Intel Xeon E5-2660 v3 CPU and 64 GB RAM.

\subsection{Non-Adaptive Interpolating PCEs - Comparison against Competitive Methods}
\label{subsec:comparison}

\begin{figure}[t!]
	\begin{tabular}[b]{c}
		\begin{tikzpicture}
		\begin{semilogyaxis}[width=0.45\textwidth, xlabel={Model evaluations}, ylabel=$\epsilon_{\text{cv}, \text{RMS}}$, legend pos=north east, ytick={10, 1, 0.1, 0.01, 0.001, 0.0001, 0.00001}, legend style={draw=none, fill=none}, title={Ishigami Function (3D)}]
		\addplot[mark=*, black, thick] table[x index=1, y index=3]{plot_data/results_comparisons/ishigami_ipce_cv_errors.txt};
		\addplot[mark=o, black, thick] table[x index=1, y index=3, col sep=comma]{plot_data/results_comparisons/ishigami_psam.txt};
		\addplot[mark=o, gray, thick] table[x index=1, y index=3, col sep=comma]{plot_data/results_comparisons/ishigami_ols2.txt};
		\addplot[mark=*, gray, thick] table[x index=1, y index=3, col sep=comma]{plot_data/results_comparisons/ishigami_ols5.txt};
		\end{semilogyaxis}
		\end{tikzpicture}
	\end{tabular}
	\hfill
	\begin{tabular}[b]{c}
		\begin{tikzpicture}
		\begin{semilogyaxis}[width=0.45\textwidth, xlabel={Model evaluations}, ylabel=$\epsilon_{\text{cv}, \text{RMS}}$, legend pos=north east,  legend style={draw=none, fill=none}, title={Cantilever Beam (4D)}]
		\addplot[mark=*, black, thick] table[x index=1, y index=3]{plot_data/results_comparisons/cantilever_ipce_cv_errors.txt};
		\addplot[mark=o, black, thick] table[x index=1, y index=3, col sep=comma]{plot_data/results_comparisons/cantilever_psam.txt};
		\addplot[mark=o, gray, thick] table[x index=1, y index=3, col sep=comma]{plot_data/results_comparisons/cantilever_ols2.txt};
		\addplot[mark=*, gray, thick] table[x index=1, y index=3, col sep=comma]{plot_data/results_comparisons/cantilever_ols5.txt};
		\end{semilogyaxis}
		\end{tikzpicture}
	\end{tabular}	
	\\
	\begin{center}
		\begin{tabular}[b]{c}
			\begin{tikzpicture}
			\begin{semilogyaxis}[width=0.45\textwidth, xlabel={Model evaluations}, ylabel=$\epsilon_{\text{cv}, \text{RMS}}$, legend pos=outer north east, title={Meromorphic Function (5D)}]
			\addplot[mark=*, black, thick] table[x index=1, y index=3]{plot_data/results_comparisons/mero5_ipce_cv_errors.txt};
			\addplot[mark=o, black, thick] table[x index=1, y index=3, col sep=comma]{plot_data/results_comparisons/mero5_psam.txt};
			\addplot[mark=o, gray, thick] table[x index=1, y index=3, col sep=comma]{plot_data/results_comparisons/mero5_ols2.txt};
			\addplot[mark=*, gray, thick] table[x index=1, y index=3, col sep=comma]{plot_data/results_comparisons/mero5_ols5.txt};
			\legend{\small{Interpolation}, \small{PSP}, \small{LSR ($C=2$)}, \small{LSR ($C=5$)}}
			\end{semilogyaxis}
			\end{tikzpicture}
		\end{tabular}
	\end{center}
	\caption{Performance of different non-intrusive, non-adaptive PCE methods in terms of the relation between approximation error and cost.}
	\label{fig:nonadaptive}
\end{figure}
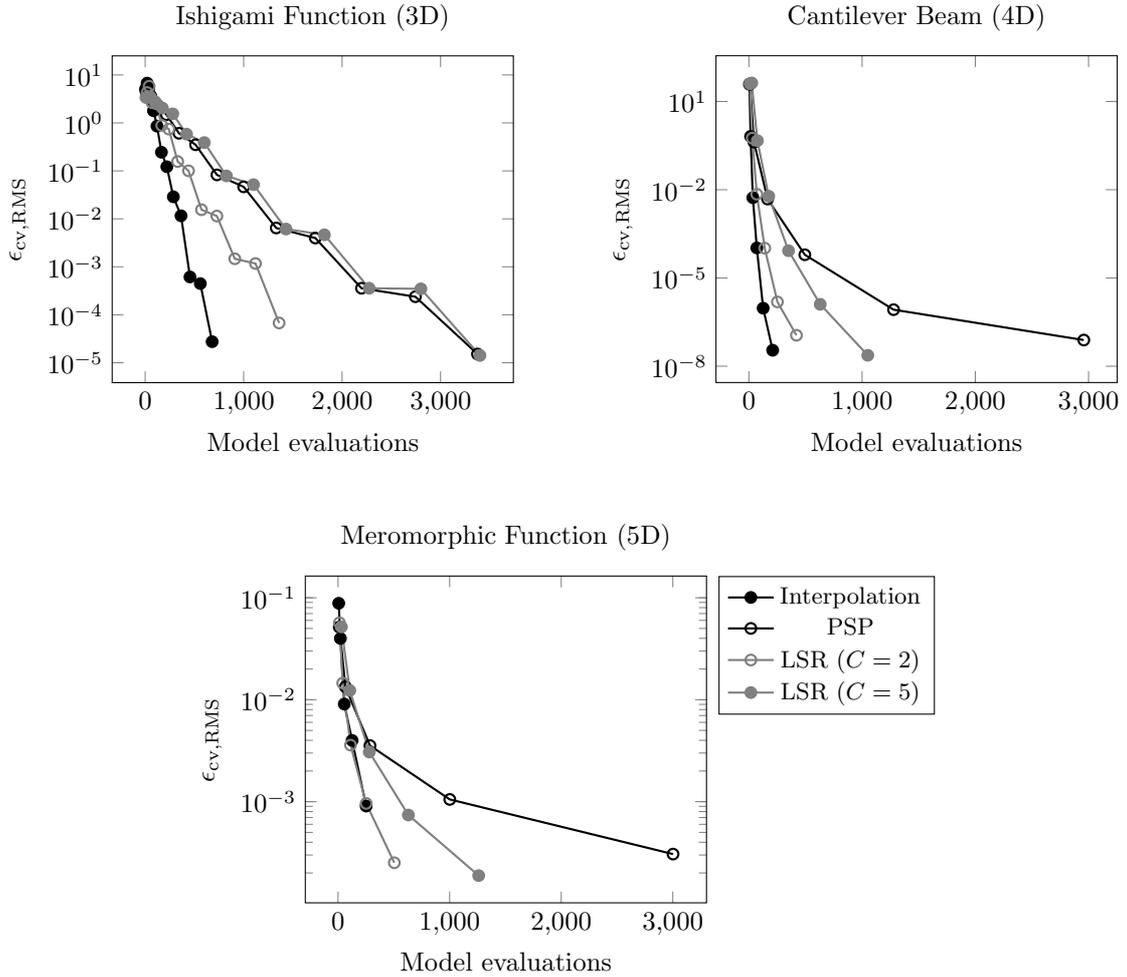

We compare interpolating \glspl{pce} against ones computed with \glsfirst{psp} and with \glsfirst{lsr}.
For both competitive methods, we employ the UQLab software \cite{marelli2014, marelli2015}.
For the \gls{lsr}-\glspl{pce}, we use two oversampling schemes, such that the experimental design is either twice or five times bigger than the polynomial basis, e.g., for an $M$-term \gls{pce}, $CM$ collocation points are used, where $C=2$ or $C=5$.
The experimental designs are based on Sobol-sequence samples, which typically result in an improved accuracy for the same number of collocation points \cite{blatman2007, loukrezisPhD}.

In this study, we do not make use of adaptivity, but instead compute \glspl{pce} based on \gls{td} bases, which correspond to the multi-index sets
\begin{equation}
\label{eq:td_set}
 \Lambda^{\mathrm{TD}}_{p_{\max}} = \left\{\mathbf{p} : \left|\mathbf{p}\right| = \sum_{n=1}^N p_n \leq p_{\max} \in \mathbb{N}_{\geq0}\right\}.
\end{equation}
In that way, for the same value of $p_{\max}$, the \glspl{pce} constructed with each method have the exact same basis polynomials. 
We investigate the performance of each \gls{pce} method in terms of the relation between approximation accuracy and cost. The approximation accuracy is measured using the \gls{rms} \gls{cv} error 
\begin{equation}
\label{eq:cverr_rms}
\epsilon_{\text{cv}, \text{RMS}} = \sqrt{\frac{1}{Q} \sum_{q=1}^Q \left(\widetilde{g}\left(\mathbf{y}_q\right) - g\left(\mathbf{y}_q\right)\right)^2}, 
\end{equation}
where $Q$ denotes the size of a \gls{cv} sample drawn randomly from the input \gls{pdf}, $\widetilde{g}$ the \gls{pce}-based surrogate model, and $g$ the original model.
In all numerical experiments we set $Q=10^5$. 
The cost refers to the number of model evaluations, equivalently, collocation points, which is necessary to compute the \gls{pce} coefficients with each approach.
Since the size of a \gls{td} basis, accordingly, the number of corresponding collocation points, grows quickly with the number of parameters, low-dimensional benchmark models are employed. 

As a first model, we use the 3-dimensional Ishigami function
\begin{equation}
\label{eq:ishigami}
g(\mathbf{y}) = \sin\left(y_1\right) + a \sin^2\left(y_2\right) + b \, y_3^4 \sin\left(y_1\right), 
\end{equation}
where $a=7$ and $b=0.1$.
The parameter vector $\mathbf{y}$ takes values in the image space of the random vector $\mathbf{Y} = \left(Y_1, Y_2, Y_3\right)$, all single \glspl{rv} of which follow the uniform distribution $\mathcal{U}\left(-\pi, \pi\right)$.

As a second model, we employ the $4$-dimensional parametric function
\begin{equation}
 \label{eq:cantilever}
 g(\mathbf{y}) = \frac{600 P_{\text{v}} + 600 P_{\text{h}}}{w t^2},
\end{equation}
which models the stress induced upon a cantilever beam of width $w$ and thickness $t$ due to the vertical and horizontal loads $P_{\text{v}}$ and $P_{\text{h}}$, respectively \cite{eldred2009}. 
The four parameters are modeled as \glspl{rv}, such that $w \sim \mathcal{N}\left(4, 10^{-4}\right)$ (in inches), $t \sim \mathcal{N}\left(2, 10^{-4}\right)$ (in inches), $P_{\text{h}} \sim \mathcal{N}\left(500, 10^4\right)$ (in N), and $P_{\text{v}} \sim \mathcal{N}\left(1000, 10^4\right)$ (in N), where $\mathcal{N}\left(\mu,\sigma^2\right)$ denotes a Gaussian distribution with mean $\mu$ and variance $\sigma^2$.

As a third model, we use the $5$-dimensional meromorphic function
\begin{equation}
\label{eq:meromorphic}
g\left(\mathbf{y}\right) = \frac{1}{1 + \mathbf{w} \cdot \mathbf{y}},
\end{equation}
where $\mathbf{w}$ is a vector of positive weights which govern the parameter anisotropy \cite{migliorati2013a}, given as $\mathbf{w} = \frac{\hat{\mathbf{w}}}{2 \|\hat{\mathbf{w}}\|_1}$, where $\hat{\mathbf{w}} = \left(1, 0.5, 0.1, 0.05, 0.001\right)$.
The parameter vector $\mathbf{y}$ takes values in the image space of the random vector $\mathbf{Y} = \left(Y_1, \dots, Y_{5}\right)$, where $Y_n \sim \mathcal{U}\left(-1, 1\right)$, $n=1,\dots,5$.

The results concerning all models are presented in Figure~\ref{fig:nonadaptive}. 
For all considered models, interpolating \glspl{pce} are clearly superior to the competitive methods in terms of their approximation accuracy for similar cost. 
The only exception is the \gls{lsr}-\gls{pce} method with an oversampling coefficient $C=2$, applied to the meromorphic function, which has a comparable performance to the interpolating \gls{pce}.

\subsection{Adaptive Interpolating PCEs}
\label{subsec:results_adaptive}

We now consider three moderately high-dimensional benchmark models, which shall be approximated using adaptive methods, specifically, with interpolating \glspl{pce} computed directly using Algorithm~\ref{algo:dimadapt_orthogonal}, hierarchical Leja interpolations based on Algorithm~\ref{algo:dimadapt_hierarchical}, as well as \glspl{pce} obtained by transforming the hierarchical Newton basis to an orthogonal one (see Section~\ref{subsec:basis_equiv}).
To depart from the standard setting of uniform or normal input \glspl{rv}, we model all inputs to follow truncated normal and Gumbel distributions. 
In that way, all methods are tested in the context of arbitrary \glspl{pdf} \cite{soize2004, wan2006}.
A truncated normal distribution shall be denoted with $\mathcal{TN}\left(\mu,\sigma^2, l, u\right)$, and is defined by a normal distribution $\mathcal{N}\left(\mu, \sigma^2\right)$ truncated in the value range $\left[l, u\right]$ \cite{burkardt2014}.
A Gumbel distribution is denoted with $\mathcal{G}\left(\ell, \beta\right)$, where $\ell$ is a location parameter and $\beta$ a scaling parameter.

As a first model, we employ the $8$-dimensional parametric function
\begin{equation}
\label{eq:borehole}
g\left(\mathbf{y}\right) = \frac{2 \pi T_{\text{u}} \left(N_{\text{u}} - N_{\text{l}}\right)}{\ln\left({\frac{r}{r_{\text{w}}}}\right) \left( 1 + \frac{T_{\text{u}}}{T_{\text{l}}} + \frac{2 L T_{\text{u}}}{\ln\left({\frac{r}{r_{\text{w}}}}\right) r_{\text{w}}^2 K_{\text{w}}}\right)},
\end{equation}
which models the water flow through a borehole (in m$^3$/yr) and is often used to test computer-based experiments, see, e.g., \cite{morris1993, narayan2014, xiong2013}.
The input parameters along with their distributions are:
\begin{enumerate}
	\item the radius of the borehole (in m), $r_\text{w} \sim \mathcal{TN}\left(0.1, 0.161812^2, 0.05, 0.15\right)$ ,
	\item the radius of influence (in m), $r \sim \mathcal{TN}\left(3700, 4900^2, 100, 50000\right)$,
	\item the transmissivity of the upper aquifer (in m$^2$/yr), $T_{\text{u}} \sim \mathcal{TN}\left(89335, 15164^2, 63070, 115600\right)$,
	\item the potentiometric head of the upper aquifer (in m), $N_{\text{u}} \sim \mathcal{TN}\left(1050, 34.64^2, 990, 1110\right)$,
	\item  the transmissivity of the lower aquifer (in m$^2$/yr), $T_{\text{l}} \sim \mathcal{TN}\left(89.5, 15.3^2, 63.1, 116\right)$,
	\item the potentiometric head of the lower aquifer (in m), $N_{\text{l}} \sim \mathcal{TN}\left(760, 34.64^2, 700, 820\right)$,
	\item the length of the borehole (in m), $L \sim \mathcal{TN}\left(1400, 161.66^2, 1120, 1680\right)$, and
	\item  the hydraulic conductivity of the borehole (in m/yr), $K_{\text{w}} \sim \mathcal{TN}\left(10950, 632.2^2, 9855, 12045\right)$.
\end{enumerate}

As a second model, we use the $10$-dimensional parametric limit state function
\begin{equation}
\label{eq:steel_column} 
g\left(\mathbf{y}\right) = F_\text{s} - P_{\text{t}} \left( \frac{1}{2 B D} + \frac{F_0 E_\text{b}}{B D H \left(E_\text{b} - P_{\text{t}}\right)}\right),
\end{equation}
which relates the reliability of a steel column to its cost \cite{eldred2009a, kuschel1997}.
The parameters and their distributions are:
\begin{enumerate}
	\item the yield stress (in MPa), $F \sim \mathcal{TN}\left(400, 1225, 295, 505\right)$,
	\item the dead weight load (in N), 	$P_\text{d} \sim \mathcal{TN}\left(500000, 25 \cdot 10^8, 350000, 650000\right)$,
	\item the first variable load (in N), $P_1 \sim \mathcal{G}\left(559495, 70173\right)$,
	\item the second variable load (in N), $P_2 \sim \mathcal{G}\left(559495, 70173\right)$,
	\item the flange breadth (in mm), 	$B \sim \mathcal{TN}\left(300, 9, 291, 309\right)$,
	\item the flange thickness (in mm), $D \sim \mathcal{TN}\left(20, 4, 14, 26\right)$,
	\item the profile height (in mm), $H \sim \mathcal{TN}\left(300, 25, 285, 315\right)$,
	\item the initial deflection (in mm), $F_0 \sim \mathcal{TN}\left(30, 100, 0, 60\right)$,
	\item Young's modulus (in MPa), $E \sim \mathcal{G}\left(208110, 3275\right)$, and
	\item the column length (in mm), $L \sim \mathcal{TN}\left(7500, 56.25, 7470, 7530\right)$.
\end{enumerate}	
In formula \eqref{eq:steel_column}, $P_{\text{t}} = P_\text{d} + P_1 + P_2$ is the total load and $E_{\text{b}} = \frac{\pi^2 E B D H^2}{2 L^2}$ is known as the Euler buckling load.

As a third model, we extend the meromorphic function \eqref{eq:meromorphic} to $16$ input parameters. 
The weight vector is modified accordingly, such that $\hat{\mathbf{w}} = \left(1, 0.5, \dots, 10^{-7}, 5 \cdot 10^{-8}\right)$.
Moreover, the input \glspl{rv} are now modeled as 
\begin{align}
\label{eq:meromorphic_dists}
Y_n \sim  
\begin{cases}
\mathcal{TN}\left(0, 1, 0, 3\right), &n=1,3,\dots,15,\\
\mathcal{TN}\left(0, 1, -3, 0\right), &n=2,4,\dots,16.
\end{cases}
\end{align}

\begin{figure}[t!]
	\begin{tabular}[b]{c}
		\begin{tikzpicture}
		\begin{semilogyaxis}[width=0.31\textwidth, xlabel=Model evaluations, ylabel=$\epsilon_{\text{cv}, \text{RMS}}$, legend pos=north east, title={Borehole (8D)}]
		\addplot[mark=o, gray, only marks] table[x index=0, y index=2]{plot_data/results_borehole_wleja/cv_results_borehole_dali.txt};
		\addplot[mark=None, gray, thick] table[x index=0, y index=2]{plot_data/results_borehole_wleja/cv_results_borehole_dali_pce.txt};
		\addplot[mark=None, black, thick] table[x index=0, y index=2]{plot_data/results_borehole_wleja/cv_results_borehole_sai_pce.txt};
		\end{semilogyaxis}
		\end{tikzpicture}
	\end{tabular}
	\begin{tabular}[b]{c}
		\begin{tikzpicture}
		\begin{semilogyaxis}[width=0.31\textwidth, xlabel=Model evaluations, legend pos=north east, title={Steel Column (10D)}]
		\addplot[mark=o, gray, only marks] table[x index=0, y index=2]{plot_data/results_steel_wleja/cv_results_steel_dali.txt};
		\addplot[mark=None, gray, thick] table[x index=0, y index=2]{plot_data/results_steel_wleja/cv_results_steel_dali_pce.txt};
		\addplot[mark=None, black, thick] table[x index=0, y index=2]{plot_data/results_steel_wleja/cv_results_steel_sai_pce.txt};
		\end{semilogyaxis}
		\end{tikzpicture}
	\end{tabular}
	\begin{tabular}[b]{c}
		\begin{tikzpicture}
		\begin{semilogyaxis}[width=0.31\textwidth, xlabel=Model evaluations, legend pos=north east, title={Meromorphic (16D)}]
		\addplot[mark=o, gray, thick, only marks] table[x index=0, y index=2]{plot_data/results_mero_wleja/cv_results_mero_dali.txt};
		\addplot[mark=None, gray, thick] table[x index=0, y index=2]{plot_data/results_mero_wleja/cv_results_mero_dali_pce.txt};
		\addplot[mark=None, black, thick] table[x index=0, y index=2]{plot_data/results_mero_wleja/cv_results_mero_sai_pce.txt};
		\end{semilogyaxis}
		\end{tikzpicture}
	\end{tabular}
	\caption{Convergence behavior of the hierarchical interpolation (gray o-marks), the basis-transformed interpolating PCE (gray solid line), and the directly computed interpolating PCE (black solid line) in terms of approximation accuracy.}
	\label{fig:cv_rms_error_convergence}
\end{figure}
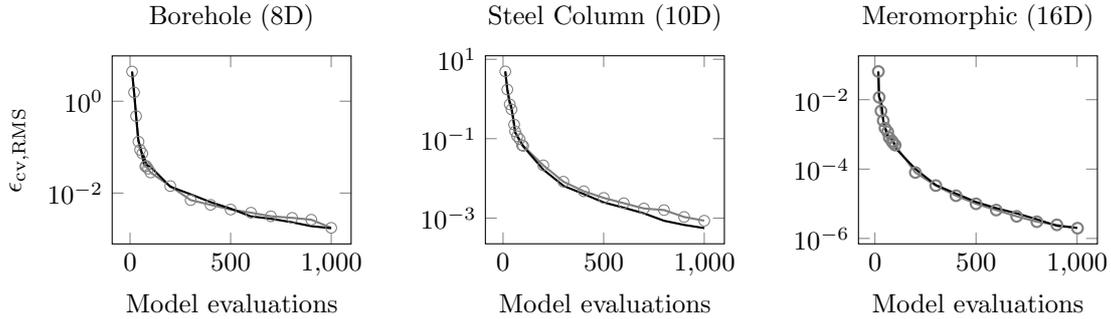

All three models are approximated by hierarchical Leja interpolations, basis-transformed interpolating \glspl{pce}, and directly computed interpolating \glspl{pce}. 
All methods are compared against one another in terms of their approximation accuracy for a similar cost, equivalently, for a similar number of collocation points.
The approximation accuracy is measured with the \gls{rms} \gls{cv} error \eqref{eq:cverr_rms}.
The results concerning all three benchmark models are presented in Figure~\ref{fig:cv_rms_error_convergence}.
As can be observed, the curves corresponding to each approach are barely distinguishable from one another. 
Therefore, all three methods can be regarded as equivalent in terms of approximation accuracy.

\begin{figure}[t!]
	\begin{tabular}[b]{c}
		\begin{tikzpicture}
		\begin{semilogyaxis}[width=0.31\textwidth, legend pos=north east, title={Borehole (8D)}, ylabel=$\epsilon_{\text{rel}, \mu}$]
		\addplot[mark=None, gray, thick] table[x index=0, y index=1]{plot_data/results_borehole_wleja/moments_relerr_borehole_dali_pce.txt};
		\addplot[mark=None, black, thick] table[x index=0, y index=1]{plot_data/results_borehole_wleja/moments_relerr_borehole_sai_pce.txt};
		\end{semilogyaxis}
		\end{tikzpicture} 
	\end{tabular}
	\begin{tabular}[b]{c}
		\begin{tikzpicture}
		\begin{semilogyaxis}[width=0.31\textwidth, legend pos=north east, title={Steel Column (10D)}]
		\addplot[mark=None, gray, thick] table[x index=0, y index=1]{plot_data/results_steel_wleja/moments_relerr_steel_dali_pce.txt};
		\addplot[mark=None, black, thick] table[x index=0, y index=1]{plot_data/results_steel_wleja/moments_relerr_steel_sai_pce.txt};
		\end{semilogyaxis}
		\end{tikzpicture} 
	\end{tabular}
	\begin{tabular}[b]{c}
		\begin{tikzpicture}
		\begin{semilogyaxis}[width=0.31\textwidth, legend pos=north east, title={Meromorphic (16D)}, ytick={0.01, 0.00001, 0.00000001}]
		\addplot[mark=None, gray, thick] table[x index=0, y index=1]{plot_data/results_mero_wleja/moments_relerr_mero_dali_pce.txt};
		\addplot[mark=None, black, thick] table[x index=0, y index=1]{plot_data/results_mero_wleja/moments_relerr_mero_sai_pce.txt};
		\end{semilogyaxis}
		\end{tikzpicture} 
	\end{tabular}
	\\
	\begin{tabular}[b]{c}
		\begin{tikzpicture}
		\begin{semilogyaxis}[width=0.31\textwidth, xlabel=Model evaluations, ylabel=$\epsilon_{\text{rel}, \sigma^2}$, legend pos=north east]
		\addplot[mark=None, gray, thick] table[x index=0, y index=2]{plot_data/results_borehole_wleja/moments_relerr_borehole_dali_pce.txt};
		\addplot[mark=None, black, thick] table[x index=0, y index=2]{plot_data/results_borehole_wleja/moments_relerr_borehole_sai_pce.txt};
		\end{semilogyaxis}
		\end{tikzpicture}
	\end{tabular}
	\begin{tabular}[b]{c}
		\begin{tikzpicture}
		\begin{semilogyaxis}[width=0.31\textwidth, xlabel=Model evaluations, legend pos=north east]
		\addplot[mark=None, gray, thick] table[x index=0, y index=2]{plot_data/results_steel_wleja/moments_relerr_steel_dali_pce.txt};
		\addplot[mark=None, black, thick] table[x index=0, y index=2]{plot_data/results_steel_wleja/moments_relerr_steel_sai_pce.txt};
		\end{semilogyaxis}
		\end{tikzpicture}
	\end{tabular}
	\begin{tabular}[b]{c}
		\begin{tikzpicture}
		\begin{semilogyaxis}[width=0.31\textwidth, xlabel=Model evaluations, legend pos=north east]
		\addplot[mark=None, gray, thick] table[x index=0, y index=2]{plot_data/results_mero_wleja/moments_relerr_mero_dali_pce.txt};
		\addplot[mark=None, black, thick] table[x index=0, y index=2]{plot_data/results_mero_wleja/moments_relerr_mero_sai_pce.txt};
		\end{semilogyaxis}
		\end{tikzpicture}
	\end{tabular}
	\caption{Convergence behavior of the basis-transformed interpolating PCE (gray solid line) and the directly computed interpolating PCE (black solid line) in terms of moment estimation accuracy.}
	\label{fig:moment_convergence}
\end{figure}
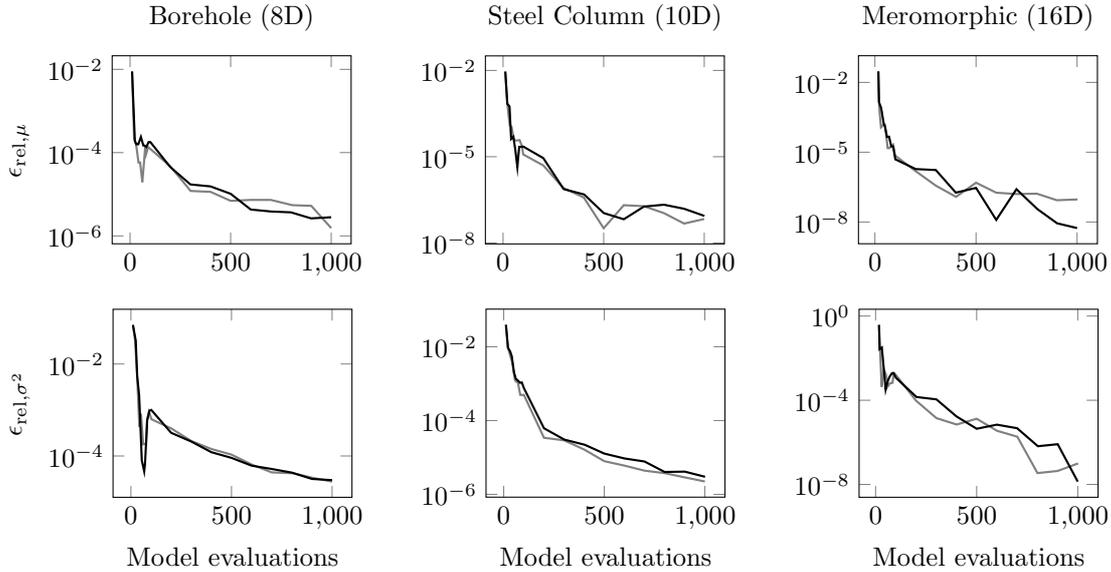

\begin{figure}[t!]
	\begin{tabular}[b]{c}
		\begin{tikzpicture}
		\begin{axis}[width=0.31\textwidth, ybar, symbolic x coords={$r_{\mathrm{w}}$, $N_{\mathrm{u}}$, $N_{\mathrm{l}}$, $L$, $K_{\mathrm{w}}$}, xtick=data, legend pos=north east, ytick={0.0, 0.2, 0.4, 0.6, 0.8}, tick label style={/pgf/number format/fixed}, bar width=5, title={Borehole (8D)}, ylabel={$S^{\text{f}}$, $S^{\text{t}}$}, xlabel=Parameters]
		\addplot[black, fill=black] coordinates {($r_{\mathrm{w}}$, 0.745)  ($N_{\mathrm{u}}$, 0.072) ($N_{\mathrm{l}}$, 0.072) ($L$, 0.069) ($K_{\mathrm{w}}$, 0.017)};
		\addplot[gray, fill=gray] coordinates {($r_{\mathrm{w}}$, 0.768)  ($N_{\mathrm{u}}$, 0.08) ($N_{\mathrm{l}}$, 0.08) ($L$, 0.077) ($K_{\mathrm{w}}$, 0.019)};
		\legend{$S^{\text{f}}_{\text{ref}}$, $S^{\text{t}}_{\text{ref}}$}
		\end{axis}
		\end{tikzpicture}
	\end{tabular}
	\hfill
	\begin{tabular}[b]{c}
		\begin{tikzpicture}
		\begin{axis}[width=0.31\textwidth, ybar, symbolic x coords={$F_{\mathrm{s}}$, $P_{\mathrm{d}}$, $P_1$, $P_2$, $D$, $F_0$}, xtick=data, legend pos=north east, ytick={0.0, 0.2, 0.4, 0.6}, tick label style={/pgf/number format/fixed}, bar width=5, title={Steel Column (10D)}, xlabel=Parameters]
		\addplot[black, fill=black] coordinates {($F_{\mathrm{s}}$, 0.624)  ($P_{\mathrm{d}}$, 0.015) ($P_1$, 0.051) ($P_2$, 0.051) ($D$, 0.186) ($F_0$, 0.068)};
		\addplot[gray, fill=gray] coordinates {($F_{\mathrm{s}}$, 0.624)  ($P_{\mathrm{d}}$, 0.015) ($P_1$, 0.052) ($P_2$, 0.052) ($D$, 0.188) ($F_0$, 0.07)};
		\end{axis}
		\end{tikzpicture}
	\end{tabular}
	\hfill
	\begin{tabular}[b]{c}
		\begin{tikzpicture}
		\begin{axis}[width=0.31\textwidth, ybar, symbolic x coords={$Y_1$, $Y_2$}, xtick=data, legend pos=north east, ytick={0.0, 0.2, 0.4, 0.6, 0.8}, tick label style={/pgf/number format/fixed}, title={Meromorphic (16D)}, xlabel=Parameters, ymax=0.8, bar width=5, enlarge x limits=0.5]
		\addplot[black, fill=black] coordinates {($Y_1$, 0.6972)  ($Y_2$, 0.2671)};
		\addplot[gray, fill=gray] coordinates {($Y_1$, 0.7212)  ($Y_2$, 0.2906)};
		\end{axis}
		\end{tikzpicture}
	\end{tabular}
	\vspace{1em}
	\\
	\begin{tabular}[b]{c}
		\begin{tikzpicture}
		\begin{semilogyaxis}[width=0.31\textwidth, ylabel=$\epsilon_{\text{rel}, S^{\text{f}}}$, legend pos=north east, title style={at={(0.8,0.8)}, anchor=north}, title={$r_{\text{w}}$}]
		\addplot[mark=None, gray, thick] table[x index=0, y index=1]{plot_data/results_borehole_wleja/sobol_f_relerr_borehole_dali_pce.txt};
		\addplot[mark=None, black, thick] table[x index=0, y index=1]{plot_data/results_borehole_wleja/sobol_f_relerr_borehole_sai_pce.txt};
		\end{semilogyaxis}
		\end{tikzpicture}
	\end{tabular}
	\hfill
	\begin{tabular}[b]{c}
		\begin{tikzpicture}
		\begin{semilogyaxis}[width=0.31\textwidth, legend pos=north east, title style={at={(0.8,0.8)}, anchor=north}, title={$F_{\text{s}}$}]
		\addplot[mark=None, gray, thick] table[x index=0, y index=1]{plot_data/results_steel_wleja/sobol_f_relerr_steel_dali_pce.txt};
		\addplot[mark=None, black, thick] table[x index=0, y index=1]{plot_data/results_steel_wleja/sobol_f_relerr_steel_sai_pce.txt};
		\end{semilogyaxis}
		\end{tikzpicture}
	\end{tabular}
	\hfill
	\begin{tabular}[b]{c}
		\begin{tikzpicture}
		\begin{semilogyaxis}[width=0.31\textwidth, legend pos=north east, title style={at={(0.8,0.8)}, anchor=north}, title={$Y_1$}]
		\addplot[mark=None, gray, thick] table[x index=0, y index=1]{plot_data/results_mero_wleja/sobol_f_relerr_mero_dali_pce.txt};
		\addplot[mark=None, black, thick] table[x index=0, y index=1]{plot_data/results_mero_wleja/sobol_f_relerr_mero_sai_pce.txt};
		\end{semilogyaxis}
		\end{tikzpicture}
	\end{tabular}
	\\
	\begin{tabular}[b]{c}
		\begin{tikzpicture}
		\begin{semilogyaxis}[width=0.31\textwidth, xlabel=Model evaluations, ylabel=$\epsilon_{\text{rel}, S^{\text{t}}}$, legend pos=north east, title style={at={(0.8,0.8)}, anchor=north}, title={$r_{\text{w}}$}]
		\addplot[mark=None, gray, thick] table[x index=0, y index=1]{plot_data/results_borehole_wleja/sobol_t_relerr_borehole_dali_pce.txt};
		\addplot[mark=None, black, thick] table[x index=0, y index=1]{plot_data/results_borehole_wleja/sobol_t_relerr_borehole_sai_pce.txt};
		\end{semilogyaxis}
		\end{tikzpicture}
	\end{tabular}
	\hfill
	\begin{tabular}[b]{c}
		\begin{tikzpicture}
		\begin{semilogyaxis}[width=0.31\textwidth, xlabel=Model evaluations, legend pos=north east, title style={at={(0.8,0.8)}, anchor=north}, title={$F_{\text{s}}$}]
		\addplot[mark=None, gray, thick] table[x index=0, y index=1]{plot_data/results_steel_wleja/sobol_t_relerr_steel_dali_pce.txt};
		\addplot[mark=None, black, thick] table[x index=0, y index=1]{plot_data/results_steel_wleja/sobol_t_relerr_steel_sai_pce.txt};
		\end{semilogyaxis}
		\end{tikzpicture}
	\end{tabular}
	\hfill
	\begin{tabular}[b]{c}
		\begin{tikzpicture}
		\begin{semilogyaxis}[width=0.31\textwidth, xlabel=Model evaluations, legend pos=north east, title style={at={(0.8,0.8)}, anchor=north}, title={$Y_1$}]
		\addplot[mark=None, gray, thick] table[x index=0, y index=1]{plot_data/results_mero_wleja/sobol_t_relerr_mero_dali_pce.txt};
		\addplot[mark=None, black, thick] table[x index=0, y index=1]{plot_data/results_mero_wleja/sobol_t_relerr_mero_sai_pce.txt};
		\end{semilogyaxis}
		\end{tikzpicture}
	\end{tabular}
	\caption{Top row: Reference Sobol indices. Middle/bottom rows: Convergence behavior of the basis-transformed interpolating PCE (gray solid line) and the directly computed interpolating PCE (black solid line) in terms of Sobol index estimation accuracy.}
	\label{fig:sobol_convergence}
\end{figure}
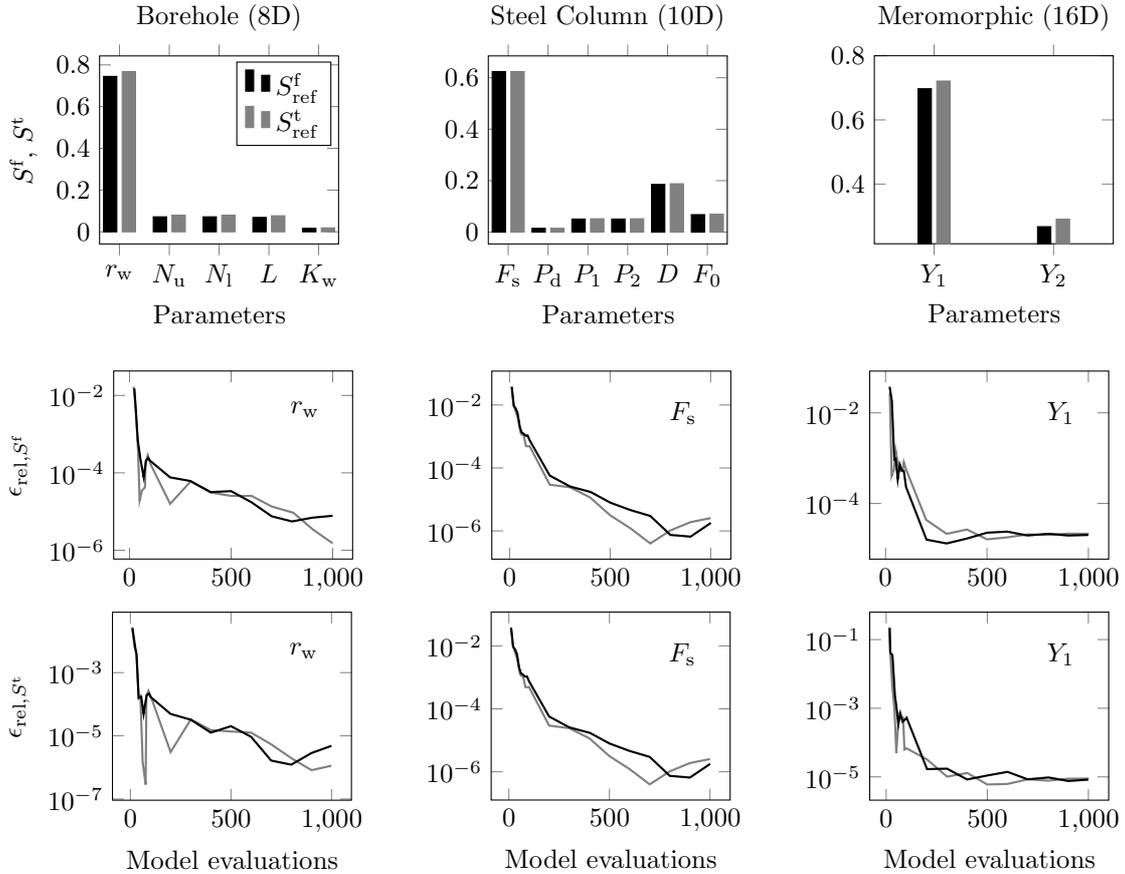

The \glspl{pce} are also post-processed for the estimation of moments and sensitivity indices.
The accuracy of the estimates is measured with a relative error metric, e.g., for the case of the expected value of the \gls{qoi},
\begin{equation}
\label{eq:relerr}
\epsilon_{\text{rel}, \mu} = \left| \frac{\mu_{\text{ref}} - \mu}{\mu_{\text{ref}}} \right|,  
\end{equation}
where $\mu_{\text{ref}}$ is a reference value, and accordingly for the other estimates.
The reference moment values are computed with quasi-\gls{mc} integration based on a Sobol-sequence sample with $10^8$ parameter realizations.
The reference Sobol index values are computed by post-processing high-order anisotropic \glspl{pce} based on a degree-adaptive \gls{lar} method \cite{blatman2011}, which is implemented in the UQLab software \cite{marelli2014, marelli2015}. 
An experimental design of $10^4$ Sobol-sequence sample points is used for the \gls{lar}-\gls{pce} method.

Figure~\ref{fig:moment_convergence} presents the performance of directly computed and basis-transformed interpolating \glspl{pce} regarding the accuracy of their moment estimates.
Minor differences can be observed depending on the model and on the moment which is estimated.
Nevertheless, the overall performance of both approaches can be regarded as comparable.
Figure~\ref{fig:sobol_convergence} shows a similar comparison, now regarding the accuracy of the Sobol index estimates.
The first row of Figure~\ref{fig:sobol_convergence} presents the reference Sobol indices for the non-negligible parameters of both models. 
We consider a parameter to be of negligible contribution if its Sobol index value is less than $0.01$.
The convergence results corresponding to the first- and total-order Sobol indices of the most important parameter of each model are exemplarily shown.
The convergence behavior regarding the remaining Sobol indices is similar.
Once more, the two approaches are comparable in terms of their estimation accuracy for an equal cost.

\section{Summary and Conclusions}
\label{sec:concl}
This work proposed a conceptually simple method to compute an interpolating \gls{pce} which features one polynomial term per collocation point.
The main component of the suggested approach is the use of Leja sequences as univariate interpolation grids.
The interpolating \gls{pce} can be derived in two ways. 
One can first construct a hierarchical interpolation based on Newton-like polynomials and then transform the polynomial basis to an orthogonal one, exploiting the one-to-one map between the two bases.
Alternatively, the interpolating \gls{pce} can be computed by directly using the orthogonal basis, without the aforementioned in-between steps.
In both cases, dedicated dimension-adaptive algorithms have been presented, shown in Algorithm~\ref{algo:dimadapt_hierarchical} and  Algorithm~\ref{algo:dimadapt_orthogonal}, respectively.

Once available, an interpolating \gls{pce} can be used for approximation, \gls{uq}, and \gls{sa} purposes.
A series of numerical experiments have verified that interpolating \glspl{pce} constitute accurate surrogate model themselves, and that their post-processing yields accurate estimations of moments and sensitivity indices. 
Comparisons between non-adaptive interpolating \glspl{pce} and ones based on \glsfirst{psp} and \glsfirst{lsr} methods show that the former typically result in a superior approximation accuracy for an equal cost.
Comparing adaptively constructed interpolating \glspl{pce} against adaptive, hierarchical, Newton-based interpolations reveals a comparable performance.
In terms of moment and sensitivity index estimation accuracy, comparisons between directly computed, adaptive interpolating \glspl{pce} and ones derived with a transform from the Newton basis to an orthogonal one show an equivalent performance. 

Due to its simplicity, the method proposed in this work can easily be combined with gradient or adjoint-based enhancement \cite{butler2013, jakeman2015},  with multi-element methods \cite{wan2006a}, or it can be applied in the context of multi-fidelity \cite{peherstorfer2018}.
Such extensions shall be pursued in later studies.

\paragraph{Acknowledgements:}
The first author would like to acknowledge the support of the German Federal Ministry for Education and Research (Bundesministerium f\"ur Bildung und Forschung - BMBF) through the research contract 05K19RDB. 
Both authors acknowledge the support of the Graduate School of Excellence Computational Engineering at the Technische Universit\"at Darmstadt. 

\bibliographystyle{siam}
\bibliography{references_aLeja}
\end{document}